\documentclass[12pt]{article}

\def\ra{\rightarrow}

\def\ss{\subseteq}

\def\e{\epsilon}
\def\d{\delta}

\def\O{\Omega}

\def\v{{\bf v}}
\def\w{{\bf w}}

\def\dist{\hbox{dist}\,}

\def\Aut{\hbox{Aut}\,}

\def\dbar{\overline{\partial}}

 \def\HollowBox #1#2{{\dimen0=#1 \advance\dimen0 by -#2       
       \dimen1=#1 \advance\dimen1 by #2                       
        \vrule height #1 depth #2 width #2                    
        \vrule height 0pt depth #2 width #1                   
        \llap{\vrule height #1 depth -\dimen0 width \dimen1}%
       \hskip -#2                                             
       \vrule height #1 depth #2 width #2}}                   
 \def\BoxOpTwo{\mathord{\HollowBox{6pt}{.4pt}}\;}             

\def\endpf{\hfill $\BoxOpTwo$}

\font\teneufm=eufm10
\font\seveneufm=eufm7
\font\fiveeufm=eufm5
\newfam\eufmfam
\textfont\eufmfam=\teneufm
\scriptfont\eufmfam=\seveneufm
\scriptscriptfont\eufmfam=\fiveeufm

\newfam\msbfam
\font\tenmsb=msbm10  scaled \magstep1 \textfont\msbfam=\tenmsb
\font\sevenmsb=msbm7 scaled \magstep1 \scriptfont\msbfam=\sevenmsb
\font\fivemsb=msbm5  scaled \magstep1 \scriptscriptfont\msbfam=\fivemsb
\def\Bbb{\fam\msbfam \tenmsb}

\def\RR{{\Bbb R}}
\def\CC{{\Bbb C}}

\newtheorem{theorem}{Theorem}[section]
\newtheorem{corollary}[theorem]{Corollary}
\newtheorem{proposition}[theorem]{Proposition}

\usepackage{graphicx}

\usepackage{amsmath}

\begin{document}

\begin{center}
{\huge \bf A New Invariant Metric}
\bigskip \\
{\huge \bf and Applications}
\bigskip \bigskip \\
Steven G. Krantz\footnote{Supported in part by a grant from
the National Science Foundation and a grant from the Dean of
Graduate Studies at Washington University.}
\end{center}
\vspace*{.25in}

\begin{quote}
{\bf Abstract:}  \sl
We study a new construction of an invariant metric for compact subgroups of 
the automorphism group of a domain in complex space.  Applications are provided.
\end{quote}

\section{Introduction}

Classically speaking, the first invariant metric in the subject
of complex analysis was the Poincar\'{e} metric (see [KRA2]) on the
disc $D$ in the complex plane $\CC$.  If one posits that the
Euclidean unit vector ${\bf e} = (1, 0) = 1 + i0$ has length
one, and applies the fact that the automorphism group (i.e, the
group of biholomorphic or conformal self-maps) of the disc acts
transitively (and transitively on directions), then one is
forced by invariance to conclude that the Poincar\'{e} metric must be
$$
F_P^D(z; \xi) = \frac{|\xi|}{1 - |z|^2} \, .
$$
[Here $\xi$ is a tangent vector and $|\xi|$ is the Euclidean length of $\xi$.]
The Poincar\'{e} metric has constant curvature $-4$ and
is a delightful and useful tool in complex function theory.
See, for instance, [KRA2].
Certainly Lars Ahlfors made good use of it [AHL] in his
study of the Schwarz lemma.

One important feature of the Poincar\'{e} metric on the disc
is that it is complete.  Put in other words, if $P \in D$ is
a fixed point and $\gamma: [0,1) \ra D$ is a $C^1$ curve
such that $\gamma(0) = P$ and $\lim_{t \ra 1^-} |\gamma(t)| = 1$
then the Poincar\'{e} length of $\gamma$ is $+\infty$.  In particular,
the metric is not smooth, nor even continuous, on $\overline{D}$; and there
is no evident way to extend the metric to $\partial D$.  One may feel
that this is part and parcel of having a transitive automorphism group
on $D$ (that is, any point can be moved to any other point by some
automorphism), and that feeling would be correct.  We shall say more
about transitive automorphism groups below.

There are a number of devices for creating invariant metrics on other
domains in $\CC$ and $\CC^n$, and also on complex manifolds.  For domains
in $\CC$, the uniformization theorem may be used to good effect---see [BER].
For this theorem gives us a covering map from $D$ to the given domain (provided
that the given domain is not $\CC$ nor $\CC \setminus \{0\}$).  And the metric
can be pushed down via the map.

Both in the plane and in higher dimensions there are constructions
of Bergman, Carath\'{e}odory, and Kobayashi/Royden.  See [KRA1], [KRA3], [KRA5].
On smoothly bounded domains in $\CC$, on strongly pseudoconvex domains
in $\CC^n$, on finite type domains in $\CC^2$, and in a variety
of other commonly encountered circumstances, these metrics are complete
(this is true regardless of whether the automorphism group acts transitively).
Thus the metric blows up at the boundary---often at the rate of the reciprocal
of the distance to the boundary.  The work [FLE] provides a relevant example
in which this last property does {\it not} hold.

It is a matter of some interest, and of effective utility as well, to construct
an invariant metric on a complex domain such that the metric is smooth
on the closure of the domain.  The purpose of this paper is to discuss and
construct such a metric, and to describe some applications.  A preliminary
version of this metric appeared in [GRK3], and we certainly recommend that
source for a profound application of the idea.

\section{Notation and Background}

We say that $\Omega \ss \CC^n$ is a {\it domain} if it is open and connected.
In this paper our domains will usually be bounded.
A mapping $\Phi: \O \ra \O$ is called an {\it automorphism} if it
is holomorphic, one-to-one, and onto.  It is automatic (but nontrivial to see)
that such a mapping automatically has a holomorphic inverse (see [KRA1]).
The collection of automorphisms of a fixed domain $\Omega$ forms a group
under the binary operation of composition of mappings.  We call this
the {\it automorphism group}, and we denote it by $\Aut(\Omega)$.

One of the features of the automorphism group that makes it
useful is the attendant topology.  We equip $\Aut(\Omega)$ with
the compact-open topology, which is equivalent to the topology
of uniform convergence on compact sets.  It is known (see [KOB1]) that,
when $\Omega$ is bounded, then $\Aut(\Omega)$ is a real Lie group.
It is never a complex Lie group.

We say that $\Aut(\Omega)$ {\it acts transitively} on $\Omega$
if, whenever $P, Q \in \O$ then there is an automorphism
$\Phi$ such that $\Phi(P) = Q$. Domains with transitive
automorphism group are relatively rare; the bounded symmetric
domains of Cartan are typical of these (see [HEL]). Among
smoothly bounded domains, the only domain with transitive
automorphism group is the unit ball (up to biholomorphic
equivalence). See the discussion in the next section as well
as [KRA1], [ROS], [WON].

A slightly weaker, and geometrically more natural, condition than
transitive automorphism group is {\it noncompact automorphism group}.
The automorphism group is noncompact if it is noncompact in
the indicated topology.  The following result of H. Cartan gives a useful
criterion for non-compactness:

\begin{proposition} \sl
Let $\Omega \ss \CC^n$ be a bounded domain.  The automorphism group
of $\O$ is noncompact if and only if there exist a $P \in \O$ and
a point $X \in \partial \O$ and a sequence of automorphisms $\varphi_j$
such that
$$
\lim_{j \ra \infty} \varphi_j(P) = X \, .
$$
\end{proposition}

If $D$ is the unit disc in $\CC$ then we may let $P = 0$, $X = 1 + i0$ and
$\varphi_j(\zeta) = (\zeta + (1 - 1/j))/(1 + (1 - 1/j)\zeta)$ to see that
$\Aut(D)$ is noncompact.  Obversely, if $A$ is the annulus then one may
see using elementary arguments that no such $P$ and $X$ exist; thus
$\Aut(A)$ is compact.  The only smoothly bounded planar domain with noncompact
automorphism group is the disc $D$ (see [KRA4]).   In higher
dimensions, the smoothly bounded domains with noncompact automorphism group
have yet to be classified.   However the only strongly pseudoconvex such
domain is the unit ball.

\section{The Original Construction on a Strongly Pseudoconvex Domain}

The original construction of a smooth-to-the-boundary invariant metric,
due to Greene and Krantz [GRK3], was on a strongly pseudoconvex domain
with smooth boundary. See [KRA1] for definitions and background on such
domains. A fundamental result for these types of domains is due to Bun
Wong [WON] and Rosay [ROS]:

\begin{theorem}  \sl
Let $\Omega \ss \CC^n$ be strongly pseudoconvex with $C^2$ boundary.
If the automorphism group of $\Omega$ is noncompact then $\Omega$ is
biholomorphic to the unit ball $B$.
\end{theorem}

We cannot prove this result here, but see [KRA1].  

Another key fact for us, that is proved in [GRK3], is the following:

\begin{theorem} \sl
Let $\Omega \ss \CC^n$ be a smoothly bounded, strongly pseudoconvex domain that is {\it not}
biholomorphic to the unit ball $B$.  Let $M$ be a positive integer.
Then there is a finite constant $K$, depending on $M$ such that
$$
\sup_{\varphi \in {\rm Aut}(\Omega), \atop
                |\alpha| \leq M, z \in \Omega} \left |
         \left ( \frac{\partial^\alpha}{\partial z^\alpha} \right ) \varphi (z) \right | \leq K \, .
$$
\end{theorem}

\noindent The proof of this result that appears in [GRK3] uses Bergman representative
coordinates, for which see [GKK].  This argument is rather elaborate, and we can only
give an indication of the idea here.
\smallskip \\

\noindent {\bf Sketch of the Proof of Theorem 3.2:}  Fix a smoothly bounded, strongly pseudoconvex
domain $\O = \{\zeta \in \CC: \rho < 0\}$ which is not biholomorphic to the ball.\footnote{We call
$\rho$ a {\it defining function} for the domain $\O$.}  If $K \ss \O$ is a compact set then
the sort of bounds described here, on the set $K$, are an immediate consequence of the Cauchy estimates.
Thus we may concentrate our attentions at the boundary of $\O$.  We now enunciate a sequence
of steps leading to Theorem 3.2.
\smallskip \\

\noindent {\bf Step 1:}  If $\e > 0$ then the number
$$
\inf\{\dist(\alpha(p), \partial \O) : \alpha \in \Aut(\O), p \in \O, \dist_{\rm Eucl}(p, \partial \O) \geq \e\}
$$
is positive.  This follows from a simple normal families argument.
\smallskip \\

\noindent {\bf Step 2:}  If $\e > 0$ then there is a $\d > 0$ such that
$$ 
\sup \bigl \{\dist_{\rm Eucl} (\alpha(p), \partial \O): \alpha \in \Aut(\O), p \in \O, \dist_{\rm Eucl}(p, \partial \O) \leq \d \bigr \} 
$$
is less than $\e$.  This also follows from a normal families argument, much like that used in Step 1.
\smallskip \\

\noindent {\bf Step 3:}  The Bergman kernel $K = K_\Omega$ for $\O$ extends smoothly to the set \\
$\overline{\O} \times \overline{\O} \setminus \{\hbox{boundary diagonal}\}$.  This is
a classical result of Kerzman [KER].  It follows from the hypoellipticity of
the $\overline{\partial}$-Neumann operator.
\smallskip \\

\noindent {\bf Step 4:}  For $w \in \partial \O$, define the {\it Levi polynomial}
$$
X(z,w) = \rho(w) + \sum_{j=1}^n (z_j - w_j) \cdot \frac{\partial \rho}{\partial w_j} \biggr |_w +
     \sum_{j,k=1}^n \frac{1}{2} (z_j - w_j)(z_k - w_k) \cdot \frac{\partial^2 \rho}{\partial w_j \partial w_k} \biggr |_w  \, .
$$
The reference [KRA1] contains detailed information about the Levi polynomial.  Then Fefferman's
asymptotic expansion for the Bergman kernel says that there are smooth functions $\varphi, \psi$ such that
$$
K_\O(z,w) = X^{-(n+1)}(z,\overline{w}) \varphi(z,w) + \bigl [ \log X(z,\overline{w}) \bigr ] \psi(z,w) \, .
$$
\vspace*{.12in}
						    
\noindent {\bf Step 5:}  Let $p \in \partial \O$.  Then there are positive numbers
$\e$ and $\eta$ such that if $z, w \in \O$ and $\dist(z,w) < \e$ and $\dist(w, p) < \e$
then $|K_\O(z,w)| \geq \eta$.   This is immediate by inspection of Fefferman's asymptotic
expansion for the Bergman kernel (Step 4).		  
\smallskip \\

\noindent {\bf Step 6:}  If $\alpha \in \Aut(\O)$ then let $J_\alpha(z)$ denote the determinant
of the complex Jacobian of $\alpha$ at the point $z \in \O$.  There is a constant $C > 0$ so that
$$
\sup \bigl \{ |J_\alpha(z)| : \alpha \in \Aut(\O), z \in \O \bigr \} < C \, ,
$$
and
$$
\inf \bigl \{ |J_\alpha(z)| : \alpha \in \Aut(\O), z \in \O \bigr \} > C^{-1} \, .
$$
This is a first step in the sort of derivative bounds that we seek.  The proof of this
result is a combination of normal families and the Cauchy estimates, using Step 5.
\smallskip \\

\noindent {\bf Step 7:}  If $p \in \O$ is a fixed point then of course $K(p,p) > 0$.  Hence
$K(z, w) \ne 0$ for $z, w$ near $p$.  So we may define, for $j = 1, \dots, n$,
$$
b_{j, p(z)} \equiv \frac{\partial}{\partial \overline{w}_j} \log \frac{K(z,w)}{K(w,w)} \biggr |_{w = p} \, .
$$
The mapping
$$
z \longmapsto (b_{1,p}(z), \dots, b_{n,p}(z))
$$
gives holomorphic coordinates near $p$.  These are the {\it Bergman representative coordinates}.  An automorphism
of $\O$, rendered in Bergman representative coordinates, will be linear.  Now there exist $\eta > 0$, $\e > 0$
such that:  If $w \in \O$ and $\dist_{\rm Eucl}(w, \partial \O) < \e$ and if $z \in \O$ with
$\dist_{\rm Eucl}(z, w) \leq \frac{3}{2} \dist_{\rm Eucl}(w, \partial \O)$, then $|K(z,w)| \geq \eta$
and $\bigl | \hbox{det}(\partial b_{j,w}/\partial z_k)_{j,k=1}^n (z)\bigr | \geq \eta$.  
The proof is a combination of Step 5 and a careful calculation with the Fefferman asymptotic
expansion (Step 4).
\smallskip \\

\noindent {\bf Step 8:}  Let $\Phi: \O_1 \ra \O_2$ be a biholomorphic mapping
of domains.  Then the Bergman kernels $K_{\Omega_1}$ and $K_{\Omega_2}$ of these domains
are related by
$$
\hbox{det\, Jac\,}_\CC \, \Phi(z) \cdot K_{\O_2}(\Phi(z), \Phi(w)) \cdot \overline{\hbox{det\, Jac\,}_\CC \,\Phi(w)} =
    K_{\O_1}(z,w) \, .
$$
\vspace*{.12in}

\noindent {\bf Step 9:} The final argument for our derivative-bound result
is by contradiction. If there is a sequence of points on which some
derivative of some automorphisms blows up, then one can use Bergman
representative coordinates (Step 7) together with the usual transformation laws for
the Bergman kernel (Step 8) to see that the Jacobians must blow up. But we have
seen in Step 6 that that is impossible. This contradiction completes the
proof. 
\endpf
\medskip \\
			  
Now, with these two results in hand, the construction of our new invariant metric
is simplicity itself.  Fix a smoothly bounded $\Omega \ss \CC^N$ that is {\it not}
biholomorphic to the unit ball.  Thus, by Theorem 3.1, the automorphism group
of $\O$ is compact.  Denote the group by $G$.  

\begin{theorem} \sl
Let $\O$ be a smoothly bounded, strongly pseudoconvex domain in $\CC^n$, and assume that
$\O$ is {\it not} biholomorphic to the unit ball in $\CC^n$.   Then there is a Riemannian
metric $h$, smooth on $\overline{\Omega}$, which is invariant under the holomorphic automorphism group
of $\Omega$.
\end{theorem}
{\bf Proof:}
Let $g$ be {\it any}
smooth Riemannian metric on $\overline{\O}$---say the Euclidean metric $g_{jk} = \delta_{jk}$. 
Define, for $1 \leq j, k \leq n$, 
$$
h({\bf v},{\bf w}) = \int_{\alpha \in  G} g(\alpha_* {\bf v}, \alpha_* {\bf w}) \, d\alpha \, ,
$$
where $d\alpha$ denotes the bi-invariant Haar measure on the compact group $G$.
Now it is immediate from our uniform estimates on derivatives (Theorem 3.2) that
$h$ is smooth on $\overline{\Omega}$ and further that $h$ is invariant under
the action of the automorphism group of $\Omega$.
\endpf 
\smallskip \\

In fact there is considerably more that can be said about this new metric $h$.  It can be arranged
(and we shall discuss the details of this assertion below) that the constructed smooth
metric on $\overline{\O}$ actually be a product metric near the boundary.  Perhaps even more
significantly, we can specify that the isometry group of $h$ be precisely
the same as the isometry group of the Bergman metric.  This result is quite useful in practice, for
it is known [KOB1] that the isometries of the Bergman metric consist of the biholomorphic self-mappings
(i.e., the automorphisms) and the conjugate biholomorphic self-mappings.

It is worth noting that there is a more general statement than Theorem 3.3 that has considerable utility
in practice:

\begin{theorem} \sl
Let $\O$ be a smoothly bounded, strongly pseudoconvex domain in $\CC^n$.  Let $G$ be a compact
subgroup of $\Aut(\Omega)$.  Then there is a Riemannian
metric $h$, smooth on $\overline{\Omega}$, which is invariant under $G$.
\end{theorem}

The proof of Theorem 3.4 is just the same as that for Theorem 3.3, so no further commentary is
needed.  But Theorem 3.4 will come up in our discussions below.

All of these devices were originally harnessed by Greene and Krantz in [GRK3] (see also
[GRK1], [GRK2]) to prove the following
result:

\begin{theorem} \sl Let $\Omega_0$ be a smoothly bounded, strongly pseudoconvex domain.
In fact write $\Omega_0 = \{z \in \CC^n: \rho_0(z) < 0\}$, where $\rho_0$ is a smooth
defining function for $\Omega_0$ (see [KRA1] for this concept).  There is a positive
integer $k$ and an $\epsilon > 0$ such that if $\Omega = \{z \in \CC^n: \rho(z) < 0\}$ and
$\|\rho - \rho_0\|_{C^k} < \epsilon$ then the automorphism group $\Aut(\Omega)$ is a subgroup
of $\Aut(\Omega_0)$.  Moreover, there is a smooth diffeomorphism ({\it not} a biholomorphism)
$\Phi: \Omega \ra \Omega_0$ such that the mapping
$$
\Aut(\Omega) \ni \varphi \longmapsto \Phi \circ \varphi \circ \Phi^{-1} \in \Aut(\Omega_0)
$$
is an injective group homomorphism.
\end{theorem}

In fact the special smooth-to-the-boundary invariant metric discussed here was used
in [GRK3] to reduce the question in Theorem 3.5 to an analogous question
for Riemannian metrics on a smooth, compact manifold.  Then a classical result
of David Ebin [EBI] could be invoked to get the desired semicontinuity.
Nowadays there are other approaches to the matter, including one based on
normal families that is due to Y. Kim [YKIM1], [YKIM2].

We cannot provide any of the details of the proof of Theorem 3.5.  The interested reader
may refer to [GRK3].  Our purpose in the present paper is to showcase the special
metric that was used, and to describe some applications.

We shall close this section by briefly outlining some of the 
previously described useful modifications of the
construction of the smooth-to-the-boundary metric $h$.	 Again we refer the reader
to [GKR3] for all the details.
\smallskip \\

\noindent {\bf The Isometry Group of \boldmath $h$ Can be Taken to Coincide with the
Isometry Group of the Bergman Metric:} \ \ Let $h$ be the smooth-to-the-boundary
invariant metric that we have constructed above.  For $\e > 0$, define
$$
\O_\e = \{z \in \O: \hbox{dist}_h(z, \partial \O) > \e\} \, .
$$
If $\e > 0$ is small then the implicit function theorem tells that $\O_\e$ is a 
strongly pseudoconvex domain
with $C^\infty$ boundary.  Each automorphism of $\O$ maps $\O_\e$ to itself
isometrically and biholomorphically.  Thus the isometry group for the Bergman
metric on all of $\O$ is isomorphic to the isometry group of the Bergman metric
of $\O$ restricted to $\O_\e$.  

Now choose $\e > 0$ so small that $\O \setminus \overline{\O_\e}$ is a tubular neighborhood
of $\partial \O$ {\it in the metric} $h$.  Let $\eta: \RR \ra \RR$ be a $C^\infty$ function
such that
\begin{enumerate}
\item $0 \leq \eta \leq 1$;
\item $\eta(x) = 1$ if $x \geq 2\e/3$;
\item $\eta(x) = 0$ if $x \leq \e/3$.
\end{enumerate}
With $b$ denoting the Bergman metric on $\O$, and $p$ any point of $\O$, we set
$$
H \biggr |_p = \bigl [ 1 - \eta(\dist_h(p, \partial \O)) \bigr ] h \biggr |_p +
                  \bigl [ \eta(\dist_h(p, \partial \O)) \bigr ] b \biggr |_p \, .
$$
Then, by inspection, the isometry group of $H$ is no larger than and no smaller
than the isometry group of the Bergman metric on $\O$.  In other words, the isometry groups
are the same.  In what follows we typically replace $h$ by $H$.
\endpf
\smallskip \\

\noindent {\bf The Metric \boldmath $H$ Can be Modified so That it is a Product Near the Boundary:} \ \ 
The metric $H$, by construction, coincides with the original smooth-to-the boundary metric $h$ near
the boundary of $\O$.  So it suffices to show that we can modify $h$ to have a product structure
near the boundary.  Let $U$ be an $h$-tubular neighborhood of $\partial \O$ and let $\pi: U \ra \partial \O$
be the $h$-projection.  Define a new metric $H^*$ on $U$ by
$$
H^*(\v, \w) = h(\pi_* \v, \pi_* \w) + (d \, \dist_h)(\v) \cdot (d \, \dist_h)(\w) \, .
$$
This new metric $H^*$ is $C^\infty$ on $\overline{U}$ because $h$ is such.  It is obviously
invariant under isometries of the Bergman metric of $\O$.  And it clearly has a product structure
on $U$.  Finally let $\delta > 0$ be small, so that the tubular neighborhood $U$ has diameter
at least $3\delta$.  Select a smooth cutoff function $\mu$ (analogous to $\eta$ above) so that
\begin{enumerate}
\item $0 \leq \mu \leq 1$;
\item $\mu(x) = 1$ if $x \geq 2\delta/3$;
\item $\mu(x) = 0$ if $x \leq \delta/3$.
\end{enumerate}
Now set
$$
\widetilde{H}\bigr |_p = \bigl [ \mu(\dist_h(p)) \bigr ] h \bigr |_p + \bigl [1 - \mu(\dist_h(p)) \bigr ] H^* \bigr |_p \, ,
$$
Then $\widetilde{H}$ is the metric that we seek---a product
near the boundary and invariant for the isometry group of the Bergman metric of $\O$.
\endpf 
\smallskip \\

The following points are worth noting, and they will play a role in what follows.
\begin{itemize}
\item The choice of $\delta$ in the last step of our construction gives a tubular boundary
neighborhood on which the metric is a product.  
\item If we choose $\epsilon = 2\delta$ in the
penultimate step, then we get a second ``layer'' just outside the product neighborhood; on that
layer, the metric comes from averaging over the automorphism group, as in Theorem 3.3.  
\item Lastly, by our construction of $H$, there is a region
interior to $\O$ on which the metric coincides with the Bergman metric for $\O$.  
\end{itemize}
If we choose $\epsilon$ and $\delta$ to be sufficiently small (with
$\epsilon = 2\delta$), then the first two layers will both lie in the
tubular neighborhood $U$ of $\partial \O$ that was specified in our
arguments.  Also an outer layer of the Bergman-metric-region will lie in $U$.
We will encounter this information again in our considerations below.

\section{Applications}

One immediate application of Theorem 3.3 is the following.

\begin{theorem} \sl
Let $\O \ss \CC^n$ be smoothly bounded and strongly pseudoconvex.  Let $\Phi: \Omega \ra \Omega$
be an automorphism.  Fix a point $P \in \partial \O$.  Suppose that $\Phi(P) = P$ and
also that
$$
d\Phi(P) = \hbox{id} \, ,
$$
where \ \ id \ \ is the usual diagonal identity matrix.  Then $\Phi(z) \equiv z$.
\end{theorem}
{\bf Proof:}  The proof is almost immediate.  We can think of the automorphism group
as acting on $\partial \O$ as a Riemannian manifold.  Since $\Phi$ fixes $P$ and
has first derivative equal to the identity, it preserves all geodesics emanating from $P$.
It follows that $\Phi$ is the identity on a relatively open subset of the boundary; hence it is
the identity on all of the boundary.  As a result, $\Phi$ is identically equal
to the identity on all of $\overline{\O}$.
\endpf
\smallskip \\

This last is a version of the famous Cartan uniqueness theorem for holomorphic mappings.
The classical proof of Cartan's result (see, for instance, [KRA1]) uses normal families
in a decisive way.  Of course normal families are not available when we are doing analysis
on the boundary.  So new methods are required. That is what we supply with our
new metric.

The following result of Bedford [BED] (proved independently by Barrett [BAR]) is now of
some interest for us.  See [KRA1] for background on the concept of finite type.

\begin{proposition} \sl
Let $\Omega \ss \CC^n$ be a smoothly bounded, pseudoconvex domain of finite type.
Let $G$ be the automorphism group of $\Omega$.  Then the action mapping
\begin{eqnarray*}	   
G \times \overline{\Omega} & \longrightarrow & \overline{\Omega} \\
(\varphi, z) & \longmapsto & \varphi(z) \\
\end{eqnarray*}
is jointly $C^\infty$.
\end{proposition}

This is an aesthetically pleasing proposition, but it also gives the uniform bound
on derivatives that we enunciated earlier for strongly pseudoconvex domains
in Theorem 3.2.  As a result, we can now conclude the following:

\begin{theorem} \sl
Let $\Omega \ss \CC^n$ be a smoothly bounded, finite type domain in $\CC^n$ with compact
automorphism group.  Then there is a Riemannian metric $h$ on $\Omega$ which is
\begin{enumerate}
\item[{\bf (i)}]  Smooth on $\overline{\O}$;
\item[{\bf (ii)}]  Invariant under the biholomorphic automorphisms of $\O$.
\end{enumerate}
\end{theorem}

As an immediate corollary we have

\begin{corollary} \sl
Let $\O \ss \CC^n$ be smoothly bounded and finite type.  Let $\Phi: \Omega \ra \Omega$
be an automorphism.  Fix a point $P \in \partial \O$.  Suppose that $\Phi(P) = P$ and
also that
$$
d\Phi(P) = \hbox{id} \, ,
$$
where \ \ id \ \ is the usual diagonal identity matrix.  Then $\Phi(z) \equiv z$.
\end{corollary}

It is worth noting explicitly what our new metric is {\it not}:
\begin{itemize}
\item  The new metric $h$ is {\it not} preserved under arbitrary holomorphic
self-mappings of $\Omega$.  And a non-biholomorphic mapping $\varphi: \O \ra \O$ will
{\it not} be distance-decreasing.
\item The new metric $h$ is {\it not} preserved under a biholomorphic mapping
$\Psi: \O_1 \ra \O_2$ of distinct domains.
\end{itemize}

We get the most information from $h$ when we study biholomorphic {\it self-maps} of a fixed
domain.	 For instance let us consider fixed points in the boundary.  We have the following result.

\begin{theorem} \sl Let $\O$ be a fixed, smoothly bounded, strongly pseudoconvex domain
in $\CC^n$ that is not biholomorphic to the ball.  
Let $P_1, \dots P_{2n}$ be points in general position in $\partial \Omega$.
If $\varphi: \O \ra \O$ is a biholomorphic self-map of $\O$ that fixes each $P_j$ then
$\varphi$ is the identity mapping.  
\end{theorem}
{\bf Proof:}  If $P_1, \dots, P_n$ are in general position (avoiding a cut locus configuration)
then they determine barycentric coordinates in the metric $h$ on $\partial \O$.  Thus $\varphi$ fixes
an entire relatively open subset in $\partial \O$, and so $\varphi \equiv \hbox{id}$.
\endpf
\smallskip \\

Some comments are in order.  It is known that it takes three fixed points in the boundary of
the unit disc $D \ss \CC$ in order to force a conformal mapping to be the identity.  But of
course the disc does {\it not} have compact automorphism group.  For the annulus (which {\it does}
have compact automorphism group), one fixed point in the boundary forces the identity.

Consider a strongly pseudoconvex domain which consists of the
unit ball $B$ in $\CC^2$ with a small, circularly symmetric bump
added in a neighborhood of $(1,0)$. It is known (see [LER])
that the only automorphisms of such a domain are rotations in
the $z_2$ variable.  A rotation in the $z_2$ variable will fix
{\it all points} of the form $(e^{i\theta},0)$ in the boundary, and
not necessarily be the identity.  But of course such points are not
in general position.

Now let $\O = \{z \in \CC^n: \rho(z) < 0\}$.  Let $\rho_0(z) = |z|^2 - 1$ be the
standard defining function for the unit ball in $\CC^n$.  Assume that
$\|\rho - \rho_0\|_{C^\infty}$ is small.   Then, by a theorem in [GK1], [GK2] the
Bergman metric of $\O$ will have negative curvature.  Also, a generic $\O$ of
this type will {\it not} be biholomorphic to the ball.  So its automorphism
group is compact, and a classical theorem of Cartan/Hadamard shows that there is a common
fixed point $P$ for the automorphism group $\Aut(\O)$. 

By our construction of the metric $\widetilde{H}$, this metric will have a layer ${\cal P}$
near the boundary on which the metric is a product metric, a layer ${\cal A}$ on which
the metric comes from averaging,
and an interior region ${\cal B}$ on which the metric coincides with the Bergman
metric.  Suppose that we have chosen $\epsilon$ and $\delta$ in the construction
of $\widetilde{H}$ to be quite small, so that Fefferman's asymptotic expansion
(Step 4 in the proof of Theorem 3.2) is valid even on the boundary
of region ${\cal B}$.  Thus the metric $\widetilde{H}$ on the region ${\cal B}$ is,
by Fefferman's work [FEF] (out near the boundary), asymptotically like the Bergman metric for the ball.
Thus a point $X$ in ${\cal B}$ out near the boundary will be metrically different from points in ${\cal P}$.
We can be sure then that such points will not be mapped by any automorphism into ${\cal P}$.  By the
maximum principle, we can conclude that {\it no point} of ${\cal B}$ will be mapped by
any automorphism into ${\cal P}$.  In fact we may look at this matter in another way:
the regions ${\cal P}$, ${\cal A}$, and ${\cal B}$ are each defined in terms of $\widetilde{H}$-metric
distance to the boundary.  So each will be preserved under automorphisms of $\O$.  

The structure described in the last paragraph enables us to see quantitatively why
Cartan's theorem (Proposition 2.1) is true:  If $\Omega$ is a smoothly bounded,
strongly pseudoconvex domain not biholomorphic to the ball and equiped with the metric $\widetilde{H}$, and if $P \in \O$,
then it is clear that a sequence of automorphisms {\it cannot} move $P$ out to the
boundary.  For if $P$ lies in ${\cal B}$ then it must stay in ${\cal B}$.  If it
lies in ${\cal A}$ then it must stay in ${\cal A}$.  And if it lies in ${\cal P}$ then
it must stay in $P$.  But in fact it must stay, under action of the automorphism group, 
in the given layer of the product structure in which $P$ lives.

We can certainly perform our construction of $\widetilde{H}$ so that the common
fixed point $P$ lies in ${\cal B}$.  It was proved in [GRK1] that Bergman
representative coordinates (see Section 3) can be then used to create an
equivariant embedding of $\O$ into $\CC^n$ so that the automorphism group
acts on $\O$ as a subgroup of the unitary group $U(n)$ acting naturally on space.
      
Let $B(P,r)$ be any $\widetilde{H}$-metric
ball centered at $P$.  Of course any $\varphi \in \Aut(\O)$ will
map $B(P,r)$ isometrically to itself.  And balls that lie entirely
in ${\cal B}$ will, by the equivariant re-embedding in the last paragraph, be essentially ``round'', with
the automorphism group acting on them by unitary rotation.  We thus have a new type of Schwarz lemma:

\begin{theorem} \sl
Let $\O$ be a smoothly bounded, strongly pseudoconvex domain in $\CC^n$ that is
$C^\infty$ sufficiently close to the unit ball.  Such a domain is generically
not biholomorphic to the ball, and we assume this property for $\Omega$.  With
notation as above, let $P \in \O$ be a common fixed point for the elements
of $\Aut(\O)$.  By the earlier discussion, we may take $P$ to lie in ${\cal B}$.
Then any $\widetilde{H}$-metric ball $B(P,r)$ is preserved
by elements of $\Aut(\O)$.  If the metric ball lies entirely in ${\cal B}$, and
if we take $\O$ to be equivariantly re-embedded, then the metric ball is acted
on via unitary rotation by elements of the automorphism group.
\end{theorem}

It is interesting to note that {\it any} smoothly bounded domain in the plane,
except for the unit disc of course, will have compact automorphism group---see [KRA4].
If we take the domain to be the annulus $A = \{\zeta \in \CC: 1 < |\zeta| < 2\}$, 
and if we concentrate on the closed subgroup $G$ of $\Aut(A)$ given by the rotations,
then the smooth-to-the-boundary invariant metrics $h$, $H$, and $\widetilde{H}$ for $G$ 
may all be taken to be the Euclidean metric.   For the full group $\Aut(A)$, the inversion
mapping $\zeta \mapsto 2/\zeta$ does {\it not} preserve Euclidean distance.  So the
three invariant metrics for all of $\Aut(\Omega)$ will be something different.
If instead $\O$ is the unit disc with the closed discs
$\overline{D}(1/2,1/10)$, $\overline{D}(i/2, 1/20)$, and $\overline{D}(-1/2, 1/30)$ removed,
then $\O$ is what we call {\it rigid}---that is, $\O$ has no automorphisms except
the identity mapping.  So again, $h$, $H$, and $\widetilde{H}$ can all be taken
to be the identity.  Finally, if $\O$ is the domain
$$
\O = D(0,1) \setminus \bigl [ \overline{D}(1/2, 1/10) \cup \overline{D}(-1/2, 1/10) 
\cup \overline{D}(i/2, 1/10) \cup \overline{D}(-i/2, 1/10) \bigr ] \, , 
$$
then the automorphism group consists of just four rotations.  These are of course
members of the orthogonal group.  So, again, we may take $h$, $H$, and $\widetilde{H}$ to
be the Euclidean metric.

Finally, consider the annulus $A' = \{\zeta \in \CC: 1/4 < |\zeta| < 1\}$.  The image
of this annulus under the mapping $\Phi: \zeta \mapsto [\zeta + 1/2]/[1 + \zeta/2]$ is
the doubly connected region $A = D(0,1) \setminus \overline{D}(10/21, 4/21)$. 
Of course the automorphism groups of $A$ and $A'$ are related by
$$
\Aut(A') \ni \varphi \longmapsto \Phi \circ \varphi \circ \Phi^{-1} \in \Aut(A) \, .
$$
Now we know that the automorphism group of $A'$ consists of rotations and the
inversion $\zeta \mapsto (1/4)/\zeta$.  At least for the subgroup consisting of the rotations,
the smooth-to-the boundary invariant metric may be taken to be the Euclidean metric.
But the automorphisms of $A$ will be nontrivial linear fractional transformations,
and these certainly do not respect Euclidean distance.  Thus $h$, $H$, and $\widetilde{H}$ for
$A$ will be new metrics, not anything familiar.

We cannot consider a bounded, infinitely connected domain in either $\CC$
or $\CC^n$, because such a domain will not have a smooth defining function
and therefore will not have smooth boundary by our definition. Any
finitely connected domain can, by a standard representation theorem (see
[KRA3]), be represented as the unit disc $D$ with finitely many smaller
closed discs removed. As soon as the number of excised discs is at least
two, it is not difficult to show that the automorphism group of the
resulting domain must be {\it finite} (see [HEI1], [HEI2]), hence the
group is certainly compact.   The theory of linear fractional transformations
is a great aid in studying this automorphism group.
							  
\section{Concluding Remarks}

It is a hallmark of modern differential geometry that a variety of geometric
analysis problems can be solved by the creation of a new metric.  As an instance,
Schoen and Yau [SCY] used this approach to prove the positive mass conjecture.

We have endeavored to show in the present paper the utility of a new invariant
metric for automorphism groups of domains.  We have indicated that the metric
can be constructed for strongly pseudoconvex domains and finite type domains.
It is not known whether there is such a metric on an arbitrary, smoothly
bounded, Levi pseudoconvex domain.

The work [GRK3] went further than the present paper in that it showed the
the construction of the invariant metrics $h$, $H$, and $\widetilde{H}$ can
be made to vary smoothly when the base domain $\O$ varies smoothly (in a suitable
topology on domains).   This in turn is based on a detailed study of
variation of the $\overline{\partial}$-Neumann problem and the Bergman kernel under
smooth variation of the domain.  Such information is useful, for instance, in the study of
semi-continuity of automorphism groups of domains (see [GRK3]).  

It is hoped that the work presented here will be a model for future investigations
in complex geometric analysis.	 Certainly the well-known list of invariant
metrics in complex analysis is by no means complete.

\newpage

\null \vspace*{.5in}

\section*{\sc References}
\vspace*{.2in}

\begin{enumerate}
\item[{\bf [AHL]}]  L. Ahlfors, An extension of Schwarz's lemma,
{\it Trans.\ Amer.\ Math.\ Soc.} 43 (1938), 359--364.

\item[{\bf [BAR]}]  D. E. Barrett, 
Regularity of the Bergman projection on domains with transverse symmetries,
{\it Math.\ Ann.} 258(1981/82), 441--446. 
				
\item[{\bf [BED]}]  E. Bedford, Action of the automorphisms of a smooth domain in $\CC^n$, 
{\it Proc.\ Amer.\ Math.\ Soc.} 93(1985), 232--234. 
					  
\item[{\bf [BER]}]  L. Bers, {\it Riemann Surfaces}, New York
University Press, New York, 1958.

\item[{\bf [EBI]}] D. Ebin, The manifold of Riemannian metrics,
{\it Proceedings of Symposia in Pure Mathematics}, vol. XV
(Global Analysis), AMS, Providence, RI, 1970, 11--40.

\item[{\bf [FLE]}]  J. E. Forn\ae ss and L. Lee, private communication.

\item[{\bf [GKK]}] R. E. Greene, K. T. Kim, and S. G. Krantz,
{\it The Geometry of Complex Domains}, Birkh\"{a}user
Publishing, Boston, MA, 2009, to appear.

\item[{\bf [GRK1]}]  R. E. Greene and S. G. Krantz, Stability
properties of the Bergman kernel and curvature properties of
bounded domains, {\it Recent Developments in Several Complex
Variables} (J. E. Forn\ae ss, ed.), Princeton University Press
(1979), 179--198.

\item[{\bf [GRK2]}] R. E. Greene and S. G. Krantz, Deformations
of complex structure, estimates for the $\dbar- $ equation,
and stability of the Bergman kernel, {\it Advances in Math.}
43(1982), 1--86.			 

\item[{\bf [GRK3]}]  R. E. Greene and S. G. Krantz, The
automorphism groups of strongly pseudoconvex domains, {\it
Math. Annalen} 261(1982), 425--446.

\item[{\bf [HEI1]}]  M. Heins, A note on a theorem of Rad\'{o} concerning
the $(1,m)$ conformal maps of a multiply-connected region into
itself, {\it Bull.\ Am.\ Math.\ Soc.} 47(1941), 128--130.

\item[{\bf [HEI2]}]  M. Heins, On the number of 1--1 directly conformal
maps which a multiply-connected plane region of finite connectivity
$p$ ($>2$) admits onto itself, {\it Bull.\ Am.\ Math.\ Soc.} 52(1946),
454--457.

\item[{\bf [HEL]}] S. Helgason, {\it Differential Geometry and
Symmetric Spaces}, Academic Press, New York, 1962.

\item[{\bf [YKIM1]}]  Y. W. Kim, Semicontinuity of compact group actions on
compact differentiable manifolds, {\it Arch. Math.} 49(1987), 450--455.

\item[{\bf [YKIM2]}] Y. W. Kim, The transformation groups and the isometry
groups, {\it Bull. Korean Math. Soc.} 26(1989), 47--52.
					       
\item[{\bf [KOB1]}] S. Kobayashi, {\it Hyperbolic Manifolds and
Holomorphic Mappings}, Dekker, New York, 1970.

\item[{\bf [KON]}] S. Kobayashi and K. Nomizu, On automorphisms of a K\"{a}hlerian
structure, {\it Nagoya Math.\ J.} 11(1957), 115--124.

\item[{\bf [KRA1]}]  S. G. Krantz, {\it Function Theory of
Several Complex Variables}, $2^{\rm nd}$ ed., American
Mathematical Society, Providence, RI, 2001.

\item[{\bf [KRA2]}] S. G. Krantz, {\it Complex Analysis: The
Geometric Viewpoint}, $2^{\rm nd}$ ed., Mathematical
Association of America, Washington, D.C., 2004.

\item[{\bf [KRA3]}]  S. G. Krantz, {\it Cornerstones of
Geometric Function Theory: Explorations in Complex Analysis},
Birkh\"{a}user Publishing, Boston, 2006.

\item[{\bf [KRA4]}]  S. G. Krantz, Characterizations of smooth
domains in $\CC$ by their biholomorphic self maps, {\it Am.
Math. Monthly} 90(1983), 555--557.

\item[{\bf [KRA5]}] S. G. Krantz, The Carath\'{e}odory and Kobayashi metrics
and applications in complex analysis, {\it Amer.\ Math.\ Monthly}
115(2008), 304--329.

\item[{\bf [LER]}] L. Lempert and L. Rubel, An independence
result in several complex variables, {\it Proc.\ Amer.\ Math.\
Soc.} 113(1991), 1055--1065.

\item[{\bf [ROS]}] J.-P. Rosay, Sur une characterization de la
boule parmi les domains de $\CC^n$ par son groupe
d'automorphismes, {\em Ann. Inst. Four. Grenoble} XXIX(1979),
91--97.

\item[{\bf [SCY]}] R. Schoen and S.-T. Yau, On the proof of the positive
mass conjecture in general relativity, {\it Comm.\ Math.\ Phys.} 65(1979),
45--76.

\item[{\bf [WON]}] B. Wong, Characterizations of the ball in
$\CC^n$ by its automorphism group, {\em Invent. Math.}
41(1977), 253--257.

\end{enumerate}

\end{document}